\magnification 1200

\def\eqnum#1{\eqno (#1)}
\def\blacksquare{\hbox{\vrule width 4pt height 4pt depth 0pt}}
\centerline {ON A CLASS OF SYMPLECTIC SIMILARITY TRANSFORMATION MATRICES
\footnote \dag {\it This work is supported by the Ministry of Sciences of the 
Republic of Macedonia within the contract ``Differential  geometrical and
topological problems and their applications'', under the Grant 08-2076/4, 
and JP Stre\v zevo, Bitola}}
\bigskip
\centerline {Jovan Stefanovski}
\par
\centerline {JP Stre\v zevo, bul. 1 Maj b.b., 7000 Bitola, Macedonia}
\medskip
\centerline {and}
\medskip
\centerline {Kostadin Tren\v cevski}
\par
\centerline {Institute of Mathematics, P.O.Box 162, 1000 Skopje, Macedonia}
\bigskip 
\bigskip
\par 
{\bf Abstract.} We present a class of symplectic matrices which transform 
by similarity given $2n\times 2n$ -dimensional matrix into Bunse-Gerstner 
form.
If the given matrix is skew-Hamiltonian, the transformation gives a solution 
of
an antisymmetric Riccati matrix equation, resulting from optimal
control of linear systems.
\par
\medskip
\noindent {\bf AMS(MOS) subject classification.} 15A24, 15A18, 35G05
\medskip
\par
{\bf 1. Introduction. }
\medskip
The paper is concerned with the class of $2n\times 2n$ dimensional
matrices:
$$
S = 
\left [ \matrix {S_{11} & S_{12} \cr
S_{21} & S_{22} \cr } \right ] 
\eqnum{1.1}
$$
where $S_{ij}$ $(i,j\in \{1,2\})$ are arbitrary real
$n\times n$ matrices. Let the matrix $J$ be defined by:
$$
J=\left[ \matrix{ 0 & I_n \cr -I_n & 0} \right] \ . 
$$
A matrix $S$ is Hamiltonian if 
$J^{\rm T} S J= -S^{\rm T}$. The Hamiltonian matrices come out from 
optimal control of linear systems [1],[2],[4],[5],[7]-[11],[17].
A matrix $S$ is skew-Hamiltonian if $J^{\rm T} S J= S^{\rm T}$. For 
example, if $S$ is a Hamiltonian matrix, 
then the matrix $S^2$ is skew-Hamiltonian. 
\par
We represent some further definitions and basic results. 
A matrix $U$ is symplectic if $ U^{\rm T} J  U=J $.
The symplectic similar transformations keep the Hamiltonian structure 
of matrices, i.e. if $S$ is a Hamiltonian matrix and $U$ is 
a symplectic matrix, then $U^{-1} S  U$ is also 
a Hamiltonian matrix. Analogously, it is easy to verify that 
the symplectic similar transformations also keep the 
structure of the skew-Hamiltonian matrices. The Hamiltonian matrices form Lie
algebra and the symplectic matrices form the corresponding Lie group.
\par
Bunse-Gerstner proved [6] that by application of an orthogonal symplectic
transformation one can reduce an arbitrary real matrix $S$ (1.1) in a form
$$
\left[\matrix{ S_{11}^\prime & S_{12}^\prime \cr 
S_{21}^\prime & S_{22}^\prime }\right]
\eqnum{1.2}
$$
such that its (2,1) block $S_{21}^\prime $ is 
upper-triangular and (1,1) block 
$S_{11}^\prime $ is upper-Hessenberg. For the transformation matrix a 
composition of finite number of Givens-like [16] and Hauseholder-like [16] 
symplectic orthogonal transformations is used.
\par
Van Loan proved [15] that if $S$ is a skew-Hamiltonian (squared Hamiltonian) 
matrix, then $S$ is a symplectic-orthogonal similar to a matrix of the form:
\par
$$
\left[\matrix{ S_{11}^\prime & S_{12}^\prime \cr \cr 
0 & S_{11}^{\prime {\rm T}}}\right]
\eqnum{1.3}
$$
where $S_{11}^\prime $ is a Hessenberg matrix and $S_{12}^\prime$ 
is a skew-symmetric matrix. 
For the transformation matrix the same composition of finite number of 
Givens-like and Hauseholder-like symplectic orthogonal transformations is used.
Besides Van-Loan's method [15],[3],[6], 
for this purpose, the methods from [4] or [13] may be applied.
\par
In [14] are applied Arnoldi-like matrices for similarity transformation of
$2n\times 2n$ -dimensional matrix $S$ into a form (1.2), such that the 
submatrix
$S_{21}^\prime =0$, under the condition that a minimal annulator of the
first column of the transformation matrix, with respect to $S$, is of degree 
$n$. Although for general matrices $S$ the last condition is connected with 
solving a nonlinear eigenvalue-eigenvector problem, for the skew-Hamiltonian
matrices $S$ these transformations result in $S_{21}^\prime =0$ also,
as do in Van Loan's work [15], because the minimal annulator of 
skew-Hamiltonian matrices is of degree $n$.
\par
In this paper we present a new class of symplectic matrices which 
transform $2n\times 2n$ -dimensional matrices into Bunse-Gerstner form (1.2),
instead of the Givens-like, Hausholder -like or
Arnoldi-like symplectic transformation matrices. 
We use the matrices of the form:
$$
U=\left [ \matrix {Q & 0 \cr Y Q & Q \cr } \right ], 
\eqnum{1.4}
$$
where $Q$ is an orthogonal matrix and $Y$ is a symmetric matrix. This
form arises naturally in solving of 
the antisymmetric Riccati matrix equation (ARME) [13]:
$$
-Y S_{12} Y + S_{22} Y - Y S_{11} + S_{21} = 0
$$
where $S$ is a skew-Hamiltonian matrix and $Y$ is an unknown symmetric
matrix, because the matrix $Y$ in (1.4) is already a solution of ARME, in case
the matrix (1.4) brings $S$ in block-triangular form. To compare our result
with [15], if a solution of ARME is required after the Van Loan's 
transformation [15] is performed,
one has to compute the matrix $U_{21} U_{11}^{-1}$, where $U_{11}$ and
$U_{21}$ are blocks of the transformation matrix. If the matrix 
$U_{11}$ is singular, it can not be done. 
\par
Other references are mentioned in the text.
\bigskip 
{\bf 2. Main results.} 
\medskip
Let be given the following $2n\times 2n$ block matrix 
$$
S = 
\left [ \matrix {S_{11} & S_{12} \cr
S_{21} & S_{22} \cr } \right ] \ , 
\eqnum{2.1}
$$
where $S_{ij}$ $(i,j\in \{1,2\})$ are arbitrary 
$n\times n$ matrices. In this section we present a finite
algorithm for similar transformation by a matrix $U$ of 
the matrix $S$, i.e.  
$$
U^{-1}  S U =
\left [ \matrix {S^\prime_{11} & S^\prime_{12} \cr
S^\prime_{21} & S^\prime_{22} \cr } \right ] = S^\prime \ ,
$$
where $S_{11}^\prime $ is an upper 
Hessenberg matrix and $S_{21}^\prime$ 
is an upper triangular matrix. 
\par
Let us consider the following set of regular $2n\times 2n$ 
matrices:
$$
G = \{ U = 
\left [ \matrix { L & 0 \cr
Y  L & L^{-{\rm T}} \cr } \right ] 
: Y^{\rm T}=Y \}.
$$
\noindent It is easy to verify that 
the set $G$ with the matrix multiplication is a
group, a subgroup of the group of the symplectic matrices. 
\par
Applying the similar transformation $U^{-1} S U$, we obtain
$$
U^{-1} S U = 
\left [ \matrix { L^{-1} & 0 \cr
-L^{\rm T}  Y & L^{\rm T} \cr } \right ]  
\left [ \matrix {S_{11} & S_{12} \cr
S_{21} & S_{22} \cr } \right ]  
\left [ \matrix { L & 0 \cr
Y  L & L^{-{\rm T}} \cr } \right ] 
$$
$$
= \left [ \matrix { L^{-1} (S_{11}+S_{12} Y) L & 
L^{-1} S_{12} L^{-{\rm T}} \cr 
L^{\rm T} (S_{21}+S_{22} Y -Y S_{11}-
Y S_{12} Y) L & 
L^{\rm T} (S_{22}-Y S_{12}) L^{-{\rm T}} \cr } \right ] 
$$
$$
= \left [ \matrix {S^\prime_{11} & S^\prime_{12} \cr
S^\prime_{21} & S^\prime_{22} \cr } \right ] = S^\prime. 
\eqnum{2.2}$$
\noindent The matrix $S$ will be transformed $n-1$ times with  respect  to  
the
matrices
$$
U_{i}=
\left [ \matrix { L_{i} & 0 \cr
Y_{i}  L_{i} & L^{-{\rm T}}_{i} \cr } \right ] \ , 
\quad Y^{\rm T}_{i}=Y_{i}\ ,\quad (1\le i\le n-1).
\eqnum{2.3}
$$
By the first transformation, the first column of the  new  matrix
$S_{21}$ will become upper 
triangular and the first column  of  the  new
matrix $S_{11}$ will become Hessenberg. Continuing  this  process,  by
the $i$-th transformation, the $i$-th column of the  new  matrix $S_{21}$
will become upper triangular 
such that the previous $i-1$ columns  of
that matrix will remain upper triangular, and the  $i$-th  column  of
the new matrix $S_{11}$ will become Hessenberg, such that previous $i-1$ 
columns of that matrix will remain Hessenberg. 
\par
From
$$
U=\left [ \matrix { L_{1} & 0 \cr
Y_{1}  L_{1} & L^{-{\rm T}}_{1} \cr } \right ]  
\left [ \matrix { L_{2} & 0 \cr
Y_{2}  L_{2} & L^{-{\rm T}}_{2} \cr } \right ] \cdots 
\left [ \matrix { L_{n-1} & 0 \cr
Y_{n-1}  L_{n-1} & L^{-{\rm T}}_{n-1} \cr } \right ] = 
\left [ \matrix { L & 0 \cr
Y  L & L^{-{\rm T}} \cr } \right ] , 
\eqnum{2.4}
$$
we obtain
$$
L = L_{1} L_{2}\cdots L_{n-1}, 
$$
$$
Y=Y_{1}+L^{-{\rm T}}_{1} Y_{2} L^{-1}_{1}+L^{-{\rm T}}_{1} L^{-{\rm T}}_{2} 
Y_{3} L^{-1}_{2} L^{-1}_{1}+\cdots 
+L^{-{\rm T}}_{1} L^{-{\rm T}}_{2}\cdots L^{-{\rm T}}_{n-2} Y_{n-1} 
L^{-1}_{n-2}\cdots L^{-1}_{2} L^{-1}_{1}
$$
\noindent and  hence $Y$  is  a  symmetric  matrix.  Further  we  will  show
explicitly the structures of $L_{i}$  and $Y_{i}$  which  are  needed  for
defining the $i$-th step of the algorithm. For the  sake  of
simplicity, the indices of the matrices $Y_i$ will be omitted.
\par
Let us put
$$
Y = \alpha  v v^{\rm T} = \alpha 
\left [ \matrix { 0 & 0 \cr
0 & \bar{v}  \bar{v}^{\rm T} \cr } \right ] \ , 
\quad v = 
\left [ \matrix {v_{1} \cr \vdots \cr v_{n} \cr } \right ] 
= \left [ \matrix {0 \cr \bar{v} \cr } \right ] \ , 
\quad \bar{v} = \left [ \matrix 
{ v_{i+1} \cr \vdots \cr v_{n} \cr } \right ], 
\eqnum{2.5}
$$
$$
L^\prime = 
\left [ \matrix {I_{i} & 0 \cr 0 & \bar{L} \cr } \right ] \ , 
\quad \bar{L} = 
\left [ \matrix { 1 & 0 & 0 & \cdots & 0 \cr 
w_{i+2} & 1 & 0 & \cdots & 0 \cr
w_{i+3} & 0 & 1 & \cdots & 0 \cr
\vdots & \vdots & \vdots & \ddots & \vdots \cr
w_{n} & 0 & 0 & \cdots & 1 \cr } \right ]\ , \quad 
w= \left [ \matrix {w_{1} \cr \vdots \cr w_{i+1} \cr 
w_{i+2} \cr \vdots \cr w_{n} \cr } \right ] = 
 \left [ \matrix {0 \cr \vdots \cr 0 \cr 
w_{i+2} \cr \vdots \cr w_{n} \cr } \right ], 
$$
where $\alpha $ is a scalar, $v$ and $w$ are vector-columns. 
If the $i$-th column of $S_{21}$ is already upper triangular,
we put $Y=0$. The matrix $Y$ is symmetric.
\par 
Let us denote 
$$
S_{11} = [r_{1},\ldots ,r_{n}] = 
\left [ \matrix {r_{11} & \cdots & r_{1n} \cr 
\vdots & \ddots & \vdots \cr 
r_{n1} & \cdots & r_{nn} \cr } \right ]\ , \quad 
S_{21} = [t_{1},\ldots ,t_{n}] = 
\left [ \matrix {t_{11} & \cdots & t_{1n} \cr 
\vdots & \ddots & \vdots \cr 
t_{n1} & \cdots & t_{nn} \cr } \right ] \ . 
$$
The $i$-th column of the matrix $S^\prime_{21}$ is upper triangular,
if and only if
$$
t_{i}-\alpha  v  v^{\rm T}  r_{i} = 
\left [ \matrix { * \cr \vdots \cr * \cr 0 \cr \vdots  \cr  0 
\cr } \right ] \ , 
$$
i.e. the last $n-i$ elements of this vector 
should be equal to zero. If 
$$
\sum ^{n}_{j=i+1} t_{ji}r_{ji} \neq 0, 
\eqnum{2.6}
$$
\noindent the  solution  of  this  system  of 
$\alpha $ and $v$ yields to 
$$
v_{j} = t_{ji} \ , \quad (i+1 \le j \le n) \qquad \hbox{    and   } 
\qquad \alpha = {1\over v^{\rm T} r_{i}} = 
{1\over \sum^{n}_{j=i+1} t_{ji} r_{ji}}. 
$$ 
\par 
The vector $w$ which takes part in the construction of the matrix 
$L^\prime$, will be determined from the condition that 
the $i$-th column of the matrix $L^{-1} S_{11} 
L$ be Hessenberg. It will be done as follows. If that 
column is already Hessenberg, we put $L=I$. 
If it is not Hessenberg, then among the elements 
$r_{i+1,i},r_{i+2,i},\ldots ,r_{ni}$ we choose that  one  of 
maximal module. Let $r_{ki}$ be such an element. For such $k$ 
we do a similar transformation of the matrix $S_{11}$ which 
consists of interchanging the $k$-th with the $(i+1)$-th row and 
also the $k$-th with the $(i+1)$-th column of $S_{11}.$ 
Thus  the  maximal  module  element  will  come  on  position 
$(i+1,i)$. For the sake of simplicity, the obtained matrix 
will also be denoted by $S_{11}$, and the matrix $L^\prime$ which 
will be obtained further should be multiplied   from  left 
with the permutation matrix in order to obtain the matrix 
$L$. The permutation matrix on $S$ belongs to the group $G$. 
\par
The $i$-th column of $L^{-1} S_{11}  L$ and 
hence of $S^\prime_{11}$ is Hessenberg if and only if:
$$
\left [ \matrix { 1 & 0 & 0 & \cdots & 0 \cr 
-w_{i+2} & 1 & 0 & \cdots & 0 \cr
-w_{i+3} & 0 & 1 & \cdots & 0 \cr
\vdots & \vdots & \vdots & \ddots & \vdots \cr
-w_{n} & 0 & 0 & \cdots & 1 \cr } \right ]  
\left [ \matrix { r_{i+1,i} \cr r_{i+2,i} \cr 
s \cr r_{ni} \cr } \right ] = 
\left [ \matrix { * \cr 0 \cr \vdots \cr 
0 \cr } \right ] \ . 
$$
Hence the solution for the vector $w$ is given by  
$$ 
w_{i+2} = {r_{i+2,i}\over r_{i+1,i}} \ , \ldots ,
w_{n} = {r_{ni}\over r_{i+1,i}}.
$$
\par 
We should prove that the first $i-1$ columns of the new obtained 
matrix $S_{21}$  remain  upper  triangular  after  the  $i$-th 
transformation, assuming that the first $i-1$ columns of the 
matrix $S_{21}$ are upper triangular. Indeed, 
$$
S^\prime _{21}=L^{\rm T} S_{21} L + \alpha  L^{\rm T}  
S_{22}  v v^{\rm T} L - \alpha  L^{\rm T}  
v v^{\rm T} L L^{-1} S_{11} L - 
\alpha ^{2}  v^{\rm T}  S_{12}  v  L^{\rm T} 
 v  v^{\rm T}  L. 
$$ 
\noindent Since the first $i$ elements of the vector row 
$v^{\rm T} L$ are equal to zero,  the  first  $i$  columns  of  the 
second and the fourth matrices  in  the  previous  sum,  are 
zeros. 
The first $i-1$ columns of the matrix $L^{-1} 
S_{11} L$ are Hessenberg  according  to  the  inductive 
assumption, and the $i$-th column has been done Hessenberg 
using the vector $w$. Thus the first $i-1$ elements of the 
following vector-row 
$$
v^{\rm T} L  L^{-1}  S_{11}  L 
$$
\noindent are zeros. Consequently, the first $i-1$ columns
of the matrix $S^\prime_{21}$ will be upper triangular. The $i$-th 
column of that matrix is already upper triangular, 
using $\alpha$ and $v$. 
\par 
We should prove that the  first  $i$  columns  of  the 
matrix  $S^\prime_{11}$  are  Hessenberg.  It  follows  from   the 
following representation 
$$
S^\prime_{11} = L^{-1} S_{11} L + 
\alpha  L^{-1}  S_{12}  v v^{\rm T} L 
$$
and from the fact that the first $i$ columns 
of the matrix $L^{-1} S_{11} L$ are Hessenberg. 
\par 
According to the following forms of the matrices $L_{i}$ 
and $Y_{i}$ 
$$
L_{1} = \left [ \matrix { 
1 & 0 & 0 & \cdots & 0 \cr 
0 & * & * & \cdots & * \cr 
0 & * & * & \cdots & * \cr 
\vdots & \vdots & \vdots & \ddots & \vdots \cr 
0 & * & * & \cdots & * \cr } \right ]\ , \quad 
L_{2} = \left [ \matrix { 
1 & 0 & 0 & \cdots & 0 \cr 
0 & 1 & 0 & \cdots & 0 \cr 
0 & 0 & * & \cdots & * \cr 
\vdots & \vdots & \vdots & \ddots & \vdots \cr 
0 & 0 & * & s & * \cr } \right ] \ , \ldots 
$$
$$
Y_{1} = \left [ \matrix { 
0 & 0 & 0 & \cdots & 0 \cr 
0 & * & * & \cdots & * \cr 
0 & * & * & \cdots & * \cr 
\vdots & \vdots & \vdots & \ddots & \vdots \cr 
0 & * & * & \cdots & * \cr } \right ] \ , \quad 
Y_{2} = \left [ \matrix { 
0 & 0 & 0 & \cdots & 0 \cr 
0 & 0 & 0 & \cdots & 0 \cr 
0 & 0 & * & \cdots & * \cr 
\vdots & \vdots & \vdots & \ddots & \vdots \cr 
0 & 0 & * & \cdots & * \cr } \right ] \ ,  \ldots 
$$
it follows that the first row and the first column 
of the matrix $Y$ are zero. 
\par
We have proved the following theorem: 
\par 
\proclaim Theorem 2.1.
Let be given a matrix $S$ 
by (2.1). Applying the similar transformations on 
$S$ by the matrices  $U_{i} ,\ i=1,\ldots ,n-1$ (2.3), 
under the conditions (2.6), we get the global transformation 
matrix $U$ (2.4) with matrix $Y$ in which the first row and 
column are equal to zero, and the transformed matrix $S^\prime$ 
be such that $S^\prime_{11}$ is an upper
Hessenberg matrix and $S^\prime_{21}$ is an  
triangular matrix. \blacksquare 
\par
{\it Remark.} One feature of the presented algorithm  is  that  it
contains only linear algebraic operations:  adding,  subtracting,
multiplying and dividing. Thus the  theorem is true also  if  the
elements of the matrix $S$ are elements of larger set. For example
it can be the field of complex numbers or the ring of the analytical 
functions of more variables.
\par
However, if the elements of the  matrix $S$  are  real  numbers,  from
numerical viewpoint it is more convenient if the matrix $L$ is
orthogonal. Therefore, we define the following algorithm.
\par
For arbitrary given matrix $S$, we shall construct the matrix: 
$$
U = \left [ \matrix {Q & 0 \cr Y Q & Q \cr } \right ] \ , 
\quad Q^{\rm T} Q=I \ , \quad Y^{\rm T}=Y
\eqnum{2.7}
$$
\noindent such that 
$$
\left [ \matrix {Q^{\rm T} & 0 \cr -Q^{\rm T} Y & Q^{\rm T} \cr } \right ]  
\left [ \matrix {S_{11} & S_{12} \cr S_{21} & S_{22} \cr } \right ]  
\left [ \matrix {Q & 0 \cr Y Q & Q \cr } \right ] = 
\left [ \matrix {S^\prime_{11} & S^\prime_{12} \cr S^\prime_{21} & 
S^\prime_{22} \cr } \right ] , 
$$
\noindent where $S^\prime_{11}$ is an upper Hessenberg matrix and 
$S^\prime_{21}$ is an upper
triangular matrix, and the first row and column of $Y$ are zeros.
\par
It verifies that the  set  of  the  matrices $U$ with 
respect to the matrix multiplication is a group. We apply 
$n-1$ such similar transformations to $S$ determined by
$$
U_{i}=
\left [ \matrix {Q_{i} & 0 \cr Y_{i} Q_{i} & Q_{i} \cr 
} \right ] \ , \quad 
Q^{\rm T}_{i} Q_{i}=I \ , \quad Y^{\rm T}_{i}=Y_{i} \ , \quad (1\le i\le n-1).
\eqnum{2.8}
$$
\noindent By the first transformation, the first column of the 
obtained matrix  $S_{21}$ will become  upper triangular,
and  the  first
column of the obtained matrix $S_{11}$ will become Hessenberg.
Generally, applying the $i$-th transformation, the $i$-th  column  of
the obtained matrix $S_{21}$ will become  
upper triangular such  that the
first $i-1$ columns will remain upper triangular and the $i$-th  column
of the obtained matrix $S_{11}$ will become Hessenberg  such  that  the
first $i-1$ columns of that matrix will remain Hessenberg. 
\par
The structures of the matrices $Y_{i}$ are the  same  as  in  the
proof of Theorem 2.1, and the matrices $Q_{i}$ 
without indices, for simplicity, are given by
$$
Q = I_{n} -2 w w^{\rm T} =
\left [ \matrix {I_{i} & 0 \cr 
0 & I-2 \bar{w}  \bar{w}^{\rm T} \cr } \right ] \ , \quad 
w=
\left [ \matrix {w_{1} \cr \vdots \cr w_{i}  \cr  w_{i+1}\cr 
\vdots \cr w_{n} \cr } \right ] = 
\left [ \matrix { 0 \cr \vdots \cr 0  \cr  w_{i+1}\cr 
\vdots \cr w_{n} \cr } \right ] \ , \quad 
\bar{w} = 
\left [ \matrix { w_{i+1}\cr 
\vdots \cr w_{n} \cr } \right ] \ , 
$$
$$
w^{\rm T} w = \bar{w}^{\rm T}  \bar{w} = 1,
\eqnum{2.9}$$
\noindent where $w$ is a vector-column whose first $i$ elements are zeros. The
matrix $Q$ is Householder's [16] and orthogonal.
\par
The vector $w$ which takes part in  the  construction  of  the
orthogonal matrix $Q$, will be determined such that the $i$-th column
of the matrix $Q^{\rm T} S_{11} Q$ to be Hessenberg, i.e.
$$
(I-2 \bar{w}  \bar{w}^{\rm T}) 
\left [ \matrix {r_{i+1,i}  \cr  r_{i+2,i}  \cr  \vdots  \cr 
r_{ni} \cr } \right ] = 
\left [ \matrix { * \cr 0 \cr \vdots \cr 
0 \cr } \right ] \ . 
\eqnum{2.10}
$$
The system (2.10) together with the orthogonality condition (2.9)
has the following solution for $w$:
$$
w_{i+1} = \beta  (r_{i+1,i}+s)
$$
$$
w_{i+2} = \beta  r_{i+2,i}
$$
$$
\vdots
$$
$$
w_{n} = \beta  r_{ni}
$$
where
$$
\beta  = [2 (s^{2}+r_{i+1,i} s)]^{-1/2} \ , \quad s=[ r^{2}_{i+1,i} 
+r^{2}_{i+2,i} +s + r^{2}_{ni}]^{1/2} .
$$
We should prove that the first $i-1$ columns of the new obtained 
matrix $S_{21}$ remain upper triangular after the  $i$-th  transformation,
assuming that the first $i-1$ columns of the matrix $S_{21}$  are  upper
triangular. We have
$$
S^\prime_{21}=Q^{\rm T} S_{21} Q + \alpha  Q^{\rm T} S_{22} v 
v^{\rm T} Q - \alpha  Q^{\rm T} v v^{\rm T} Q Q^{\rm T} 
S_{11} Q -\alpha ^{2} v^{\rm T} S_{12} v Q^{\rm T} 
v v^{\rm T} Q. 
$$
\noindent Since the first $i$ elements of the vector-row $v^{\rm T} Q$ are  
zeros, it
follows that the first $i$ columns of the second and 
the fourth matrices in the previous sum are zeros. The first 
$i-1$ 
columns of the matrix $Q^{\rm T} S_{11} Q$ are 
Hessenberg according to the 
inductive assumption, and the $i$-th column was made Hessenberg  using
the vector $\bar{w}$. Thus the first $i-1$ elements of the following 
vector-row 
$$
v^{\rm T} Q Q^{\rm T} S_{11} Q
$$
\noindent are zeros. Consequently, the first $i-1$ columns of the matrix 
$S^\prime_{21}$
will be also upper triangular. The $i$-th column of that  matrix  was
done upper triangular using $\alpha $ and $v$.
\par
The first row and column of the matrix $Y$ are zero, from  the
same reasons as in the proof of Theorem 2.1. 
\par
We have proved the following theorem. 
\par
\proclaim Theorem 2.2.
Let be given a matrix $S$
with (2.1). Applying the similar transformations on 
$S$ by the matrices $U_{i} ,\ i=1,\ldots ,n-1$ (2.8), 
under the conditions of finiteness of $\beta$ in each step,  
we get the global transformation 
matrix $U$ (2.7) with matrix $Y$ in which the first row and 
column are equal to zero, and the transformed matrix $S^\prime$ 
be such that $S^\prime_{11}$ is upper 
Hessenberg matrix and $S^\prime_{21}$ is upper 
triangular matrix. \blacksquare 
\par
As a consequence from these results follow  two  corollaries
for the Hamiltonian and skew-Hamiltonian matrices.
\par
\proclaim Corollary 2.3.
Let the matrix $S$ be a Hamiltonian  
matrix. Applying $n-1$ similar transformations
(2.3) of the matrix 
$S$ under the conditions (2.6), 
we get a matrix $S^\prime$,  
such that
the left lower block of dimension 
$n\times n$ 
is a diagonal matrix and the upper left block of dimension 
$n\times n$ is
of  Hessenberg  form.
\par
{\it Proof}.  It  follows  from  (2.2)  that  the  similar  matrix
$S^\prime=U^{-1} S U$ also has the 
Hamiltonian form, if $U$ is any  matrix 
of
the form given in Theorems 2.1 and 2.2. 
The rest of the proof
is obvious. \blacksquare 
\par
The second corollary is about the ARME.
\par
\proclaim Corollary 2.4.
Let $S$ be a skew-Hamiltonian matrix. Applying $n-1$ similar 
transformations
(2.3), under the condition (2.6), the matrix $S$ 
can be  reduced
in a form such that lower left block of dimension 
$n\times n$ is zero,
the upper left block of dimension 
$n\times n$ has a Hessenberg form.
If the final similar transformation is given by the matrix $U$, i.e. 
$$
U=\left[\matrix{ L & 0\cr YL & L^{-{\rm T}} }\right] \quad {\rm or}\quad
U=\left[\matrix{ Q & 0\cr YQ & Q }\right] \ , 
$$
then the  matrix $Y$ 
is a symmetric solution of the ARME
$$
-Y S_{12} Y + S_{22} Y - Y S_{11} + S_{21} = 0,
$$
such  that  its  first  row  and  column  are 
zeros.
\par
{\it Proof.} Using the similar 
transformation $S^\prime =U^{-1} S U$, where 
$U$
is any of the matrices given in Theorems  2.1  and  2.2,  
it is obtained skew-Hamiltonian matrix $S^\prime $.
The rest of the proof is obvious. \blacksquare 
\par
{\it Remark.} Breaking down the algorithm happens when there 
is no solution of ARME or there is no solution of ARME 
with first row and column equal to zero. 
\par 
{\it Remark.} If the algorithm breaks  down  or  if  it  is 
requested a solution of ARME which is non-singular matrix, 
or the matrix with given first row and column, we could introduce a new
unknown symmetric matrix $X$ by $Y=M+N^{-{\rm T}}XN^{-1}$, where $M$ is a given 
symmetric matrix and $N$ is a given non-singular matrix.
The new equation of $X$ 
is also of the ARME type: 
$$
-X S^\prime_{12} X-X S^\prime_{11}+ S^\prime_{22} X+S^\prime_{21}=0 \ , 
$$
where  
$$  
S^\prime_{11}=N^{-1} (S_{11}+S_{12}M) N=S^{\prime {\rm T}}_{22} \ , 
$$
$$ 
S^\prime_{12}=N^{-1} S_{12} N^{-{\rm T}}=-S^{\prime {\rm T}}_{12} \ , 
$$
$$ 
S^\prime_{21}=N^{\rm T} (S_{21}+S_{22}M 
-MS_{11}-MS_{12}M)  N=-S^{\prime {\rm T}}_{21} \ , 
$$
$$
S^\prime_{22}=N^{\rm T} (S_{22}-MS_{12}) N^{-{\rm T}}=S^{\prime {\rm T}}_{11}
$$
(compare with (2.2).).
\par
{\it Example.} In this example it is $n=6$ and the 
matrix $S$ is skew-Hamiltonian with randomly  
chosen numbers from $(0,1)$ and 
symmetrically from $(-1,0)$.  The  beginning
matrix $S$ is
$$
S_{11}= \left[ \matrix { 
0.001251& 0.563585& 0.193304 &0.808740 &0.585009 &0.479873\cr 
0.350291& 0.895962& 0.822840 &0.746605 &0.174108 &0.858943\cr 
0.710501& 0.513535& 0.303995 &0.014985 &0.091403 &0.364452\cr 
0.147313& 0.165899& 0.988525 &0.445692 &0.119083 &0.004669\cr 
0.008911& 0.377880& 0.531663 &0.571184 &0.601764 &0.607166\cr 
0.166234& 0.663045& 0.450789 &0.352123 &0.057039 &0.607685\cr 
} \right ] \ \ , \ \ S_{22}=S_{11}^{\rm T}
$$
$$
S_{12}= \left[ \matrix { 
0.000000 &0.783319 &0.802606& 0.519883& 0.301950 &0.875973\cr 
-0.783319& 0.000000& 0.726676& 0.955901& 0.925718 &0.539354\cr 
-0.802606& -0.726676& 0.000000& 0.142338& 0.462081& 0.235328\cr 
-0.519883& -0.955901& -0.142338& 0.000000& 0.862239& 0.209601\cr 
-0.301950& -0.925718& -0.462081& -0.862239& 0.000000& 0.779656\cr 
-0.875973& -0.539354& -0.235328& -0.209601& -0.779656& 0.000000\cr 
} \right ]
$$
$$
S_{21}=\left[ \matrix { 
0.000000& 0.843654& 0.996796& 0.999695 &0.611499 &0.392438\cr 
-0.843654& 0.000000& 0.266213& 0.297281 &0.840144 &0.023743\cr 
-0.996796& -0.266213& 0.000000& 0.375866 &0.092624 &0.677206\cr 
-0.999695& -0.297281& -0.375866& 0.000000 &0.056215& 0.008789\cr 
-0.611499& -0.840144& -0.092624& -0.056215& 0.000000& 0.918790\cr 
-0.392438& -0.023743& -0.677206& -0.008789 &-0.918790& 0.000000\cr 
} \right ]
$$
After the similar transformation (2.2), the matrix $S$ becomes
$$
S_{11}=\left[ \matrix { 
0.001252& 2.394542& 0.628127& -1.333712 &1.170301& -2.988671\cr 
-0.822757& 0.705668& -0.064203 &-0.380534 &0.078576 &-2.116390\cr 
0.000000 &4.483943 &1.443984 &3.651134 &-1.359573 &-5.058558\cr 
0.000000 &-0.000000 &-0.232568 &-1.192501 &0.923134 &-0.061005\cr 
-0.000000& -0.000000& 0.000000 &-2.882063 &1.947789 &-0.289124\cr 
-0.000000& -0.000000& -0.000000& 0.000000 &-0.112915 &-0.049842\cr 
} \right ]
$$
$$
S_{22}=S_{11}^{\rm T}
$$
$$
S_{12}=\left[ \matrix { 
0.000000& -1.299940& -0.502376& -0.391241& 0.293239 &0.453059\cr 
1.299940& -0.000000& 1.315821 &-0.067147 &0.335201 &0.359848\cr 
0.502376& -1.315821& 0.000000 &0.003709 &-0.814256 &1.134890\cr 
0.391241& 0.067147 &-0.003709 &-0.000000 &0.579926 &-0.053570\cr 
-0.293239& -0.335201 &0.814256& -0.579926 &-0.000000 &-0.078375\cr 
-0.453059& -0.359848& -1.134890 & 0.053570 & 0.078375 & 0.000000\cr 
} \right ]
$$
$$
S_{21}=0
$$
\par
The symmetric solution of  ARME  with  zero  first  row  and 
column is
$$
Y=\left[ \matrix { 
-0.000000& 0.000000 &-0.000000& 0.000000 &0.000000 &-0.000000\cr 
0.000000 &-2.976399 &0.633676 &-1.401061 &-1.558562& -0.186413\cr 
-0.000000& 0.633676 &-1.373721& -0.563505& 0.769165& -1.002085\cr 
0.000000 &-1.401061 &-0.563505& -0.599371& -1.051505& -0.065432\cr 
0.000000 &-1.558562 &0.769165 &-1.051505 &1.478251 &-2.829218\cr 
-0.000000& -0.186413 &-1.002085& -0.065432& -2.829218& 2.524725\cr 
} \right ].
$$
\medskip 
\medskip 
{\bf Conclusions.} We have presented a class of symplectic matrices which
we have used in transformation of matrices $S$ into Bunse-Gerstner form. Our
numerical experiments with skew-Hamiltonian matrices $S$ showed that for 
larger $n$ ($n>10$), the elements of the matrix $S^\prime _{21}$ are not
zeros, as theoretically they have to be. Further research could be to 
find the reasons for that behavior and to find an improvement of the
presented algorithm for annihilating the (2,1) block of skew-Hamiltonian
matrices. For example, one could try with multiple application of the
algorithm.
\medskip 
\medskip 
\centerline{ R E F E R E N C E S} 
\medskip 
\item{[1]} G. Ammar and V. Mehrmann, On Hamiltonian and symplectic 
Hessenberg forms, {\it Lin. Alg. App.}, {\bf 149} (1991), 55-72.
\item{[2]} B.D.O. Anderson and J.B. Moore, {\it Optimal control: linear 
quadratic method}, Prentice Hall, Englewood Cliffs (N.Y.USA), 1992.
\item{[3]} P. Benner, R. Byers and E. Barth, HAMEV and SQRED: Fortran 77
subroutines for computing the eigenvalues of Hamiltonian matrices using
Van Loan's square reduced method, Preprint-Reihe dos Chemnitzer SFB393,
SFB 393/96-06, May 1996.
\item{[4]} P. Benner, V. Mehrmann and H. Xu, A numerically stable, 
structure preserving method for computing the eigenvalues of real Hamiltonian 
or symplectic pencils, {\it Numer.Math.}, {\bf 78} (1998) Issue 3, 329-358.
\item{[5]} P. Benner, V. Mehrmann V. and H. Xu, A new method for computing 
the stable invariant subspace of a real Hamiltonian matrix, {\it Journal of 
Computational and Applied Mathematics}, {\bf 86} (1997), 17-43. 
\item{[6]} A. Bunse-Gerstner, Matrix factorization for symplectic $QR$-like 
method, {\it Lin.Alg.App.}, {\bf 83} (1986), 49-77. 
\item{[7]} W.A. Coppel, Matrix quadratic equations, {\it Bill. Austral. 
Math.Soc.}, {\bf 10} (1974), 377-401. 
\item{[8]} F.R. Gantmaher, {\it Theory of matrices}, Moscow, "Nauka", 1988 
(in Russian).
\item{[9]} A.J. Laub, A Schur method for solving algebraic Riccati 
equations, {\it IEEE Trans. Autom. Contr.}, AC-24, (1979), 913-925. 
\item{[10]} A.J. Laub, Invariant subspace methods for the numerical 
solution 
of Riccati equation, {\it The Riccati equation}, Eds. S. Bittanti, A.J. Laub 
and J.C. Willems, Springer-Verlag, Berlin, (1991), 163-196. 
\item{[11]} C. Paige and C.F. Van Loan, A Schur decomposition for 
Hamiltonian
matrices, {\it Lin. Alg. App.}, {\bf 41} (1981), 11-32. 
\item{[12]} J.E. Potter, Matrix quadratic solutions, 
{\it J.SIAM Appl.Math.}, {\bf 14}, No.3, May 1966. 
\item{[13]} J. Stefanovski and K. Tren\v cevski, Antisymmetric Riccati 
matrix equation, Proc. Ist Congress of the mathematicians and informaticians
of Macedonia, October 3-5, 1996, Ohrid, 83-92. 
\item{[14]} J. Stefanovski, Generating equations approach for quadratic
matrix equations, {\it Numerical Linear Algebra with Applications}, {\bf 6} 
(1999), 295-326. 
\item{[15]} C.F. Van Loan, A symplectic method for approximating all the 
eigenvalues of a Hamiltonian matrix, {\it Lin. Alg. App.}, {\bf 16}  
(1984), 233-251. 
\item{[16]} J.H. Wilkinson, {\it The algebraic eigenvalue problem}, London:
Oxford University Press, 1965. 
\item{[17]} H. Xu and L. Lu, Properties of a quadratic matrix equation 
and the solution of the continuous-time algebraic Riccati equation,
{\it Lin.Alg.App.}, {\bf 222} (1995), 127-146. 
\end